\documentclass[12pt,a4paper]{amsart}
\usepackage[top=1.15in, bottom=1.15in, left=1.17in, right=1.17in]{geometry}
\usepackage{amssymb,amsmath}
\usepackage{pdflscape}
\usepackage{amscd}
\usepackage{fancyhdr}
\usepackage{tikz}
\usetikzlibrary{matrix,arrows,cd}
\usepackage{verbatim}
\usepackage{todonotes}

\usepackage{dsfont}
\usepackage{hyperref}

\numberwithin{equation}{section}
\makeatletter
\def\section{\def\@secnumfont{\mdseries}\@startsection{section}{1}%
  \z@{.7\linespacing\@plus\linespacing}{.5\linespacing}%
  {\normalfont\scshape\centering}}
\def\subsection{\def\@secnumfont{\bfseries}\@startsection{subsection}{2}%
  {\parindent}{.5\linespacing\@plus.7\linespacing}{-.5em}%
  {\normalfont\bfseries}}
\makeatother

\newtheorem{theorem}{Theorem}

\newtheorem{lemma}[theorem]{Lemma}
\newtheorem{proposition}[theorem]{Proposition}

\theoremstyle{definition}

\allowdisplaybreaks
\title{A piecewise-linear formula for rigid representations of Dynkin quivers}
\author{Deniz Kus, Markus Reineke}

\address{ 
Ruhr-University Bochum\\ Faculty of Mathematics\\
Universit\"ats\-strasse 150\\
44780 Bochum (Germany)}
\email{deniz.kus@rub.de, markus.reineke@rub.de}

\numberwithin{equation}{section}

\begin{document}

\begin{abstract} We derive a piecewise-linear formula for the rigid representation of a Dynkin quiver of a given dimension vector, and illustrate the formula in several examples. 
\end{abstract}
\maketitle
\parindent0pt

\section{Introduction}
It is a fundamental fact in the representation theory of Dynkin quivers that any dimension vector admits  a unique (up to isomorphism) rigid representation, that is, one with vanishing group of self-extensions. This fact already plays a decisive role in proofs of Gabriel's theorem stating that the isomorphism classes of representations of a Dynkin quiver are parametrized by non-negative integer valued functions on the set of positive roots of the corresponding simply-laced root system (we recall basic notation, facts and techniques from the representation theory of quivers in Sections \ref{notation}, \ref{sub22}).\\[1ex]
It is thus natural to ask how the rigid representation of a given dimension vector is given in this parametrization. Our main result is a closed piecewise-linear formula for this (Theorem \ref{main}), which is proved using standard methods of Auslander-Reiten theory, and of the theory of general subrepresentations, in Section \ref{proof}. We illustrate the efficiency of this formula by making it explicit, and comparing it to natural rank conditions, for representations of type $A$ quivers in Section \ref{appla}.\\[2ex]
{\bf Acknowledgments:} This work originated in a research seminar on flat families of quiver Grassmannians at Ruhr-University Bochum in 2022. The authors would like to thank all participants, in particular M.~Bertozzi, M.~Boos, V.~Genz, W.~Gnedin and V.~Rappel, for their contributions.

\section{Piecewise-linear formula: the main result}
\subsection{}\label{notation} For all techniques and more details on the representation theory of quivers, we refer to \cite{Schiffler}. Let $Q$ be a finite quiver with set of vertices $Q_0$ and set of arrows $Q_1$, arrows being denoted by $\alpha:i\rightarrow j$ for $i,j\in Q_0$. Let $t,h:Q_1\rightarrow Q_0$ be the maps associating to an arrow $\alpha$ its tail and head, respectively. The associated Euler form is given by
$$\langle{\bf d},{\bf e}\rangle=\sum_{i\in Q_0}d_ie_i-\sum_{\alpha:i\rightarrow j}d_ie_j,\ \ {\bf d}=(d_i)_i,\ {\bf e}=(e_i)_i\in\mathbb{Z}Q_0.$$
For an algebraically closed field $\mathbb{K}$, we consider finite-dimensional $\mathbb{K}$-representations $V$ of $Q$,  consisting of $\mathbb{K}$-vector spaces $V_i$ for $i\in Q_0$ and $\mathbb{K}$-linear maps $V_\alpha:V_i\rightarrow V_j$ for all arrows $\alpha:i\rightarrow j$. They naturally form a $\mathbb{K}$-linear abelian category ${\rm rep}_{\mathbb{K}}Q$, which is equivalent to the category  of finite dimensional left modules over the path algebra $\mathbb{K}Q$ of $Q$. For a path $\beta$ in the path algebra we denote by $V_{\beta}$ the composition of the linear maps $V_{\alpha}$, $\alpha\in Q_1$ along the path $\beta$. The tuple ${\rm\bf dim}(V)=(\dim_{\mathbb{K}}V_i)_{i\in Q_0}\in\mathbb{N}Q_0$ is called the dimension vector of $V$. We have
$$\dim_{\mathbb{K}}{\rm Hom}_Q(V,W)-\dim_{\mathbb{K}}{\rm Ext}^1_Q(V,W)=\langle{\rm\bf dim}(V),{\rm\bf dim}(W)\rangle$$
for all representations $V,W$ of $Q$. The set 
$$R_{\bf d}(Q)=\bigoplus_{\alpha:i\rightarrow j}{\rm Hom}_{\mathbb{K}}(V_i,V_j),\ \ {\bf d}\in\mathbb{N}Q_0$$ of all representations of $Q$ on the vector spaces $V_i$ of dimension $d_i$ is an affine space, on which the algebraic group $$G_{\bf d}=\prod_{i\in Q_0}{\rm GL}(V_i)$$ acts via base change
$$(g_i)_i\cdot(V_\alpha)_\alpha=(g_jV_\alpha g_i^{-1})_{\alpha:i\rightarrow j},$$
so that the orbits $\mathcal{O}_V$ of $G_{\bf d}$ in $R_{\bf d}(Q)$ correspond bijectively to the isomorphism classes $[V]$ of $\mathbb{K}$-representations of $Q$ of dimension vector ${\bf d}$. We have
\begin{equation}\label{bri}\dim R_{\bf d}(Q)-\dim \mathcal{O}_V=\dim_{\mathbb{K}}{\rm Ext}^1_Q(V,V)\end{equation}
$$\dim G_{\bf d}-\dim R_{\bf d}(Q)=\langle{\bf d},{\bf d}\rangle.$$
\subsection{}\label{sub22} We assume in the rest of the paper that $Q$ is a Dynkin quiver (that is, the underlying unoriented graph is a disjoint union of Dynkin diagrams of type $A_n$, $D_n$, $E_6$, $E_7$, $E_8$). Then the quadratic form associated to the Euler form $\langle-,-\rangle$ is positive definite. We denote by $$\Phi^+=\{{\alpha}\in\mathbb{N}Q_0\, |\, \langle{\alpha},{\alpha}\rangle=1\}$$
the set of positive roots of this form, or, equivalently, of the root system for the Dynkin type of $Q$. By the classical theorem of Gabriel, the isomorphism classes of indecomposable $\mathbb{K}$-representations of $Q$ are in bijection to $\Phi^+$. More precisely, for every ${\alpha}\in \Phi^+$ there exists a unique (up to isomorphism) indecomposable representations $U_{\alpha}$ such that ${\rm\bf dim}(U_{\alpha})={\alpha}$. By the theorem of Krull-Remak-Schmidt, every finite-dimensional representation $V$ of $Q$ is then isomorphic to a direct sum
\begin{equation}\label{dec} V\cong \bigoplus_{{\alpha}\in\Phi^+}U_{\alpha}^{{\bf m}_V({\alpha})}\end{equation}
for a unique function ${\bf m}_V:\Phi^+\rightarrow\mathbb{N}$. In particular, for a given dimension vector ${\bf d}$, there are only finitely many isomorphism classes of $\mathbb{K}$-representations of $Q$ of dimension vector ${\bf d}$. Equivalently, there are only finitely many orbits of the connected algebraic group $G_{\bf d}$ in the irreducible algebraic variety $R_{\bf d}(Q)$, thus there exists a unique open orbit $\mathcal{O}_{\bf d,\rm{ri}}\subseteq R_{\bf d}(Q)$, corresponding to a so-called rigid representation $V_{{\bf d},\rm{ri}}$ uniquely defined up to isomorphism. It is already uniquely determined by the properties ${\rm\bf dim}(V_{{\bf d},\rm{ri}})={\bf d}$ and ${\rm Ext}^1_Q(V_{{\bf d},\rm{ri}},V_{{\bf d},\rm{ri}})=0$ (see equation \eqref{bri}). However, a formula for the decomposition \eqref{dec} of $V_{{\bf d},\text{ri}}$ into indecomposables is still not present in the literature. In special cases, the description is straightforward; for example, we have $V_{\alpha,\rm{ri}}\simeq U_{\alpha}$ for all positive roots $\alpha\in \Phi^+$.\\[1ex]
Our aim is to derive a closed, piecewise-linear, formula for the decomposition of $V_{{\bf d},\rm{ri}}$ and illustrate the formula in type $A$ leading also to maximizing certain rank tuples. We write ${\bf e}\hookrightarrow {\bf d}$ (resp.~${\bf d}\twoheadrightarrow{\bf e}$) if every $\mathbb{K}$-representation of $Q$ of dimension vector ${\bf d}$ admits a subrepresentation (resp.~factor representation) of dimension vector ${\bf e}$. In the following, we use the convenient notation $[x]_+=\max\{x,0\}$ for $x\in\mathbb{R}$. The main result is the following.

\begin{theorem}\label{main} Let $Q$ be a Dynkin quiver. For the rigid representation $V_{{\bf d},\rm{ri}}$ of $Q$ of dimension vector ${\bf d}$, we have
$$\mathbf{m}_{V_{{\bf d},\rm{ri}}}(\alpha)=\left[\min\{\langle{\bf e},{\bf d}\rangle,\langle{\bf d},{\bf e'}\rangle: 0\not={\bf e}\hookrightarrow {\alpha},\ {\alpha}\twoheadrightarrow{\bf e'}\not=0) \}\right]_+\mbox{ for all }\alpha\in \Phi^+.$$
%$$G_{\bf d}\simeq\bigoplus_{{\alpha}\in\Phi^+}U_{\alpha}^{[\min(\langle{\bf e},{\bf d}\rangle\, :\, 0\not={\bf e}\hookrightarrow {\alpha},\, \langle{\bf d},{\bf e}\rangle\, :\, {\alpha}\twoheadrightarrow{\bf e}\not=0)]_+}.$$
\end{theorem}
Before proving the theorem, we recall a few concepts from the representation theory of quivers, namely almost split sequences, and generic subrepresentations/extensions.\\[1ex]
By Auslander-Reiten theory, for every indecomposable representation $U_{\alpha}$, there exists a so-called almost split sequence
$$U_{\alpha}\stackrel{f}{\longrightarrow} X\stackrel{g}{\longrightarrow}V$$
inducing the following exact sequence of $\bmod_\mathbb{K}$-valued $\mathbb{K}$-linear functors on ${\rm rep}_{\mathbb{K}}Q$:
$$0\longrightarrow{\rm Hom}_Q(V,-)\longrightarrow{\rm Hom}_Q(X,-)\longrightarrow{\rm Hom}_Q(U_{\alpha},-)\longrightarrow S_{U_{\alpha}}\longrightarrow 0,$$
where $S_{U_{\alpha}}$ is the simple functor whose value at any representation $Y$ is the vector space $\mathbb{K}^{\mathbf{m}_Y(\alpha)}$ (the above sequence thus providing a projective resolution of this functor). In particular, for every representation $Y$, we have
$$\mathbf{m}_Y(\alpha)=\dim_{\mathbb{K}}{\rm Hom}_Q(U_{\alpha},Y)-\dim_{\mathbb{K}}{\rm Hom}_Q(X,Y)+\dim_{\mathbb{K}}{\rm Hom}_Q(V,Y).$$
Moreover, the representation $X$ fulfills ${\rm Ext}_Q^1(X,X)=0$ and is therefore rigid of dimension vector $\alpha+\mathbf{dim}(V)$. If $U_{\alpha}$ is an injective indecomposable, then $X=U_{\alpha}/{\rm soc}(U_{\alpha})$ and $V=0$. Otherwise, the almost split sequence is short exact, the representation $V$ is again indecomposable, and from the above exact sequence of functors, we see immediately that every non-split mono starting in $U_{\alpha}$ factors over $f$; dually, every non-split epi to $V$ factors over $g$.\\[1ex]
The association of $V$ to $U_{\alpha}$ is functorial, in fact given by the inverse Auslander-Reiten translation functor which, in the Dynkin case, is given as $\tau^{-1}={\rm Ext}^1_Q((\mathbb{K}Q)^*,-)$; the Auslander-Reiten formulas then imply ${\rm Hom}_Q(V,U_{\alpha})=0$. The functor $\tau^{-1}$ induces a bijection from isomorphism classes of non-injective indecomposables to isomorphism classes of non-projective indecomposables, and it induces a linear map, again denoted by $\tau^{-1}\in{\rm GL}(\mathbb{Z}Q_0)$, on the level of dimension vectors. This map satisfies
\begin{equation}\label{claim0}\langle \tau^{-1}{\bf d},{\bf e}\rangle=-\langle{\bf e},{\bf d}\rangle.\end{equation}

\subsection{}\label{proof} Given two $\mathbb{K}$-representations $V_1$ and $V_2$ of $Q$, there exists a unique (up to isomorphism) representation $Y$ where the function $\dim_{\mathbb{K}}{\rm End}(Y)$ obtains its minimum among all extensions of $V_1$ by $V_2$; we call it the generic extension of $V_1$ by $V_2$ and denote it by $V_1*V_2$. We say that a representation $W$ degenerates to a representation $W'$ if the orbit $\mathcal{O}_{W'}$ is contained in the Zariski-closure of $\mathcal{O}_W$. In particular, $V_{\bf d,\rm{ri}}$ degenerates to $W$ for any representation $W$ of dimension vector ${\bf d}$.\\[1ex]
We then have the following result from \cite[Proposition 2.4]{ReGenExt}: if $V_1$ degenerates to $V_1'$ and $V_2$ degenerates to $V_2'$, and $Y$ is an arbitrary extension of $V_1'$ by $V_2'$, then $V_1*V_2$ degenerates to $Y$.
%{\color{red} Let $V_1\leq_{\rm{deg}} V_1'$ and $V_2\leq_{\rm{deg}} V_2'$, then $V_1*V_2\leq_{\rm{deg}} Y$ for all $Y\in {\rm Ext}^1_Q(V_1',V_2')$.}
From this we conclude the following. \begin{equation}\label{if22}\text{If $V_{\bf d,\rm{ri}}$ admits a factor of dimension vector ${\bf e}$, then ${\bf d}\twoheadrightarrow{\bf e}$.}\end{equation} Indeed, if $V_{\bf d,\rm{ri}}/W'\cong W$ for some representation $W$ of dimension vector ${\bf e}$, then $V_{\bf e,\rm{ri}}$ degenerates to $W$. Thus, by the above discussion, $V_{\bf e,\rm{ri}}*W'$ degenerates to $V_{\bf d,\rm{ri}}$. But this implies $V_{\bf d,\rm{ri}}\simeq V_{\bf e,\rm{ri}}*W'$ and hence  $V_{\bf e,\rm{ri}}$ is a factor of $V_{\bf d,\rm{ri}}$. Similarly, we obtain the statement for subrepresentations.
\begin{equation}\label{if22i}\text{If $V_{\bf d,\rm{ri}}$ admits a subrepresentation of dimension vector ${\bf e}$, then ${\bf e}\hookrightarrow{\bf d}$.}\end{equation}
We will use these two statements in the rest of this paper. We can now prove the theorem.\\[1ex]
\textit{Proof of Theorem~\ref{main}:} The function $(V,W)\mapsto\dim_{\mathbb{K}}{\rm Hom}_Q(V,W)$ is upper semicontinuous on the product $R_{\bf e}(Q)\times R_{\bf d}(Q)$ and its minimal value, denoted by ${\rm hom}({\bf e},{\bf d})$, is obtained on the product of open orbits, that is, ${\rm hom}({\bf e},{\bf d})=\dim_{\mathbb{K}}{\rm Hom}_Q(V_{\bf e,\rm{ri}},V_{\bf d,\rm{ri}})$. In particular, for every $\alpha\in \Phi^+$, we have
$$\dim_{\mathbb{K}}{\rm Hom}_Q(U_{\alpha},V_{\bf d,\rm{ri}})={\rm hom}({\alpha},{\bf d}).$$ 
Moreover, using the Euler form and a formula of Schofield \cite[Theorem 5.4]{Scho}, we get
$${\rm hom}({\bf e},{\bf d})=\max\{\langle{\bf f},{\bf d}\rangle\, :\, {\bf e}\twoheadrightarrow{\bf f}\}.$$
Let $U_\alpha\stackrel{f}{\rightarrow} X\stackrel{g}{\rightarrow} V$ the almost split sequence starting in $U_\alpha$ for some $\alpha\in \Phi^+$ and let $\gamma$ and $\beta$ be the dimension vectors of $X$ and $V$, respectively. From the discussion in Subsection~\ref{sub22} we have 
$$\mathbf{m}_{V_{\bf d,\rm{ri}}}(\alpha)=\dim_{\mathbb{K}}{\rm Hom}_Q(U_{\alpha},V_{\bf d,\rm{ri}})-\dim_{\mathbb{K}}{\rm Hom}_Q(X,V_{\bf d,\rm{ri}})+\dim_{\mathbb{K}}{\rm Hom}_Q(V,V_{\bf d,\rm{ri}})$$
$$={\rm hom}({\alpha},{\bf d})-{\rm hom}(\gamma,{\bf d})+{\rm hom}(\beta,{\bf d})$$
$$=\max\{\langle{\bf e},{\bf d}\rangle\, :\, {\alpha}\twoheadrightarrow{\bf e}\}-\max\{\langle{\bf e},{\bf d}\rangle\, :\, \gamma\twoheadrightarrow{\bf e}\}+\max\{\langle{\bf e},{\bf d}\rangle\, :\, \beta\twoheadrightarrow{\bf e}\}.$$
For ${\bf e}\in\mathbb{N}Q_0$, we denote by $Q_{\bf e}$ the set of all dimension vectors ${\bf f}$ such that ${\bf e}\twoheadrightarrow{\bf f}$. We claim that
\begin{equation}\label{claim}Q_{\gamma}=(Q_{\alpha}+Q_{\beta})\setminus\{{\alpha}\}.\end{equation}
In case of $U_{\alpha}$ beeing injective, we have $\beta=0$ and $\gamma={\rm\bf dim}(U_{\alpha}/{\rm soc}(U_{\alpha}))$. So identity (\ref{claim}) follows from the fact that ${\rm soc}(U_{\alpha})$ is simple, and thus the proper factors of $U_{\alpha}$ are precisely the factors of $U_{\alpha}/{\rm soc}(U_{\alpha})$. So we assume that $U_{\alpha}$ is non-injective, and prove identity (\ref{claim}) in three steps.\\[1ex]
\textit{Step 1:} We first prove $Q_{\gamma}\subseteq (Q_{\alpha}+Q_{\beta})\setminus\{\alpha\}$. Suppose that ${\bf f}\in Q_{\gamma}$ and let $Y=X/C$ be a factor of $X$ of dimension vector ${\bf f}$. Then $U_{\alpha}/(U_{\alpha}\cap C)$ is a factor of $U_{\alpha}$, say of dimension vector ${\bf f}_1$, and $$X/(U_{\alpha}+C)\simeq(X/U_{\alpha})/((U_{\alpha}+C)/U_{\alpha})\simeq V/((U_{\alpha}+C)/U_{\alpha})$$ is a factor of $V$, say of dimension vector ${\bf f}_2$, so that ${\bf f}={\bf f}_1+{\bf f}_2$. Hence the inclusion is obtained from \eqref{if22} if we can exclude the case $\bf f =\alpha$. However, if ${\bf f}={\alpha}$ we have that $U_{\alpha}$ is already a factor of $X$, which is impossible by the above exact sequence of functors characterizing almost split sequences, namely ${\rm Hom}_Q(X,U_{\alpha})={\rm Hom}_Q(V,U_{\alpha})=0$ since ${\rm End}_Q(U_{\alpha})\simeq \mathbb{K}$.\\[1ex]
\textit{Step 2:} Next, we prove $(Q_{\alpha}\setminus\{{\alpha}\})+Q_{\beta}\subseteq Q_{\gamma}$. So let $p:U_{\alpha}\rightarrow Y$ and $q:V\rightarrow Z$ be surjective maps, where $p$ is not an isomorphism. Hence, $p$ factors through $f:U_{\alpha}\rightarrow X$ to a map $\hat{p}:X\rightarrow Y$ which is still surjective. We can also compose $q$ with the map $g:X\rightarrow V$. It is then easily verified that the map $(\hat{p},q\circ g):X\rightarrow Y\oplus Z$ is again surjective. \\[1ex]
\textit{Step 3:} Finally, we prove ${\alpha}+(Q_{\beta}\setminus\{0\})\subseteq Q_{\gamma}$. So let $Z\neq 0$ and $V\rightarrow Z$ be surjective with kernel $K\rightarrow V$ being non-split epi, thus factoring to an embedding $K\rightarrow X$, and $X\rightarrow X/K$ is the desired factor of dimension vector ${\alpha}+{\rm\bf dim}(Z)$.\\[1ex]
So identity (\ref{claim}) is proven, and we rewrite it in the following ways:
\begin{equation}\label{claim2}Q_\alpha+Q_\beta=Q_\gamma\cup\{\alpha\},\end{equation}
\begin{equation}\label{claim3}Q_\gamma=((Q_\alpha\setminus\{\alpha\})+Q_\beta)\cup(\alpha+(Q_\beta\setminus\{0\}),\end{equation}
and we will make use of the following simple identity:
\begin{equation}\label{claim4}[\min\{x,y,x+y\}]_+=[\min\{x,y\}]_+.\end{equation}
Given a subset $B\subseteq \mathbb{N}Q_0$, we abbreviate 
$$\max\langle B,{\bf d}\rangle=\max\{\langle{\bf f},{\bf d}\rangle\, :\,f\in B\}$$
and calculate:
\begin{align*}&\mathbf{m}_{V_{\bf d,\rm{ri}}}(\alpha)=\max\langle Q_\alpha+Q_\beta,{\bf d}\rangle-\max\langle Q_\gamma,{\bf d}\rangle&\\&\stackrel{(\ref{claim2})}{=}\max\{\max\langle Q_\gamma,{\bf d}\rangle,\langle\alpha,{\bf d}\rangle\}-\max\langle Q_\gamma,{\bf d}\rangle=[\min\{\langle\alpha-{\bf e},{\bf d}\rangle\, :\, \gamma\twoheadrightarrow{\bf e}\}]_+
&\\&\stackrel{(\ref{claim3})}{=}\big[\min\big\{\min\{\langle\alpha-{\bf e}-{\bf f},{\bf d}\rangle\, :\, {\bf e}\in Q_\alpha\setminus\{\alpha\},\, {\bf f}\in Q_\beta\},&\\& \hspace{7cm}\min\{\langle-{\bf f},{\bf d}\rangle\, :\, {\bf f}\in Q_\beta\setminus\{0\}\}\big\}\big]_+&\\&=\big[\min\big\{\min\{\langle{\bf e}-{\bf f},{\bf d}\rangle\, :\, 0\not={\bf e}\hookrightarrow\alpha,\; \beta\twoheadrightarrow{\bf f}\},\min\{\langle-{\bf f},{\bf d}\rangle\, :\, \beta\twoheadrightarrow{\bf f}\not=0\}\big\}\big]_+&\\&=\big[\min\big\{\min\{\langle{\bf e}-{\bf f},{\bf d}\rangle\, :\, 0\not={\bf e}\hookrightarrow\alpha,\; \beta\twoheadrightarrow{\bf f}\neq 0\},\min\{\langle{\bf e},{\bf d}\rangle\, :\, 0\not={\bf e}\hookrightarrow\alpha\},&\\& \hspace{7cm}\min\{\langle-{\bf f},{\bf d}\rangle\, :\, \beta\twoheadrightarrow{\bf f}\not=0\}\big\}\big]_+&\\&
\stackrel{(\ref{claim4})}{=}\big[\min\big\{\min\{\langle{\bf e},{\bf d}\rangle\, :\, 0\not={\bf e}\hookrightarrow\alpha\},\min\{-\langle{\bf f},{\bf d}\rangle\, :\, \beta\twoheadrightarrow{\bf f}\not=0\}\big\}\big]_+.\end{align*}
%$$=[\min(\min(\langle{\bf e}-{\bf f},{\bf d}\rangle\, :\, 0\not={\bf e}\hookrightarrow\alpha,\; \beta\twoheadrightarrow{\bf f}\not=0),$$
%$$\min(\langle{\bf e},{\bf d}\rangle\, :\, 0\not={\bf e}\hookrightarrow\alpha),\min(\langle-{\bf f},{\bf d}\rangle\, :\, \beta\twoheadrightarrow{\bf f}\not=0))]_+=$$
Denoting by $Q_\alpha^{\rm inj}$ the set of dimension vectors of injective factor representations of $U_{\alpha}$, the properties of the Auslander-Reiten translation imply that $$Q_{\beta}\setminus\{0\}=\tau^{-1}(Q_{\alpha}\setminus Q_\alpha^{\rm inj}).$$
Together with the formula (\ref{claim0}), this yields
$$\min\{-\langle{\bf f},{\bf d}\rangle\, :\, \beta\twoheadrightarrow{\bf f}\not=0\}=\min\{\langle{\bf d},{\bf e'}\rangle\, :\, {\bf e'}\in Q_\alpha\setminus Q_\alpha^{\rm inj}\},$$
and thus the calculations above give 
$$\mathbf{m}_{V_{\bf d,\rm{ri}}}(\alpha)=\big[\min\{\langle{\bf e},{\bf d}\rangle,\ \langle{\bf d},{\bf e'}\rangle :\, 0\not={\bf e}\hookrightarrow\alpha,\ {\bf e'}\in Q_\alpha\setminus Q_\alpha^{\rm inj}\}\big]_+.$$
We claim that this already equals
$$[\min\{\langle{\bf e},{\bf d}\rangle,\ \langle{\bf d},{\bf e'}\rangle :\, 0\not={\bf e}\hookrightarrow\alpha,\ \alpha\twoheadrightarrow{\bf e'}\not=0\}]_+,$$
which then proves the theorem.\\[1ex]
Indeed, suppose that $\alpha\twoheadrightarrow{\bf e'}\not=0$ for ${\bf e'}$ the dimension vector of an injective representation $I$, thus in particular ${\rm Hom}_Q(U_{\alpha},I)\not=0$. The inverse Nakayama functor $\nu^{-1}={\rm Hom}_Q((\mathbb{K}Q)^*,-)$ then yields a projective representation $\nu^{-1}I$ such that 
$${\rm Hom}_Q(U_{\alpha},I)\simeq{\rm Hom}_Q(\nu^{-1}I,U_{\alpha})^*,$$
 thus we can choose a non-zero map $0\not=\varphi:\nu^{-1}I\rightarrow U_{\alpha}$. We denote its (projective) kernel by $P$, its image by $U'$, and by ${\bf f}$ the dimension vector of $U'$ which satisfies $0\neq{\bf f}\hookrightarrow \alpha$ by \eqref{if22i}. Then
 $$0\rightarrow P\rightarrow \nu^{-1}I\rightarrow U'\rightarrow 0$$
 is a projective resolution, and we find for any dimension vector ${\bf d}$ that
 $$\langle{\bf f},{\bf d}\rangle=\langle{\rm\bf dim}(\nu^{-1}I),{\bf d}\rangle-\langle{\rm\bf dim}(P),{\bf d}\rangle=\langle{\bf d},{\rm\bf dim}(I)\rangle-\langle{\rm\bf dim}(P),{\bf d}\rangle=$$
 $$=\langle{\bf d},{\bf e'}\rangle-\langle{\rm\bf dim}(P),{\bf d}\rangle\leq\langle{\bf d},{\bf e'}\rangle,$$
 which shows that the above minimum does not change when adding the terms $\langle{\bf d},{\bf e'}\rangle$ for $0\not={\bf e'}\in Q_\alpha^{\rm inj}$.
 
\subsection{} To make our main result more efficient, we have to determine, for a given positive root $\alpha\in \Phi^+$, the dimension vectors ${\bf e}\not=0$ satisfying ${\bf e}\hookrightarrow\alpha$ and $\alpha\twoheadrightarrow{\bf e}$ respectively. 

\begin{lemma} Let $\alpha\in \Phi^+$ and ${\bf e}\not=0$ be a non-zero element of $\mathbb{N}Q_0$.  We have ${\bf e}\hookrightarrow\alpha$ if and only if 
$$[\langle\beta,\alpha\rangle]_+\geq\langle\beta,{\bf e}\rangle,\ \ \forall \beta\in\Phi^+.$$  Dually, we have $\alpha\twoheadrightarrow{\bf e}$ if and only if
$$[\langle\alpha,\beta\rangle]_+\geq\langle{\bf e},\beta\rangle\ \ \forall \beta\in\Phi^+.$$

\proof By \cite[Theorem 4.2]{KR}, we have ${\bf e}\hookrightarrow\alpha$ if and only if $\dim_{\mathbb{K}}{\rm Hom}_Q(U_\beta,U_\alpha)\geq\langle\beta,{\bf e}\rangle$ for all $\beta\in\Phi^+$. Now Dynkin quivers are representation-directed, that is, there exists a partial ordering $\preceq$ on $\Phi^+$ such that ${\rm Hom}_Q(U_\beta,U_\alpha)\not=0$ already implies $\beta\preceq\alpha$, and ${\rm Ext}^1_Q(U_\beta,U_\alpha)\not=0$ implies $\alpha\prec\beta$. This immediately gives $\dim_{\mathbb{K}}{\rm Hom}_Q(U_\beta,U_\alpha)=[\langle\beta,\alpha\rangle]_+$. The dual statement is similar and the lemma is proved.
\endproof
\end{lemma}
\section{Applications to type A: rank tuples}\label{appla}
Let $Q$ be an arbitrary type $A_n$ quiver. The positive roots are parametrized by intervals in this quiver,
$$\Phi^+=\{\alpha_{ij}\, :\, 1\leq i\leq j\leq n\},$$
where $(\alpha_{ij})_k=1$ precisely if $i\leq k\leq j$. 
We denote by $U_{i,j}$ the corresponding indecomposable representations. For a subset $I\subseteq\{1,\ldots,n\}$, we denote by ${\bf e}_I$ the dimension vector such that $({\bf e}_I)_i=1$  when $i\in I$, and zero otherwise. Then $0\not={\bf e}_I\hookrightarrow \alpha_{ij}$ if and only if $\emptyset\not=I\subseteq\{i,\dots,j\}$ is closed under successors, that is, if $k\in I$ and there is an arrow $k\rightarrow l$ in $ Q$, then $l\in I$, for $i\leq k,l\leq j$. Dually, $\alpha_{ij}\twoheadrightarrow {\bf e}_I\not=0$ if and only if $\emptyset\not=I\subseteq\{i,\dots,j\}$ is closed under predecessors.
\subsection{}
We consider the equioriented type $A_n$ quiver
$$1\rightarrow 2\rightarrow\cdots\rightarrow n.$$
From the discussion above, the only non-zero dimension vectors of subrepresentations (resp.~quotient representations) of the corresponding indecomposable $U_{i,j}$ are the $\alpha_{kj}$ (resp.~$\alpha_{ik}$) for $i\leq k\leq j$. We thus find with our main result that
$$\mathbf{m}_{V_{{\bf d},\rm{ri}}}(\alpha_{ij})=\left[\mathrm{min}\left\{d_k-d_{i-1},d_k-d_{j+1}\, :\, i\leq k\leq j\right\}\right]_+,\ \ \forall 1\leq i\leq j\leq n$$
%$$G_{\bf d}=\bigoplus_{i\leq j}U_{ij}^{[\min(d_k-d_{i-1},d_k-d_{j+1}\, :\, i\leq k\leq j)]_+},$$
where we understand $d_0=0=d_{n+1}$. The representation $V=V_{{\bf d},\rm{ri}}$, given as
$$V_1\stackrel{f_1}{\rightarrow}V_2\stackrel{f_2}\rightarrow\cdots\stackrel{f_{n-1}}{\rightarrow}V_n$$ is determined up to isomorphism by the ranks
$$r_{ij}(V)={\rm rank}(f_{j-1}\circ\cdots\circ f_i:V_i\rightarrow V_j)$$
for $1\leq i\leq j\leq n$ (where we view the empty composition of maps as the identity map, that is, $r_{ii}(V)=\dim_{\mathbb{K}}V_i$). We have $r_{ij}(V)=\min\{d_i,\ldots,d_j\}$. To see this, let $i<j$ (the case $i=j$ is true by convention) and observe
\begin{align*}r_{ij}(V)&=\dim_{\mathbb{K}}V_i-\dim_{\mathbb{K}}{\rm Ker}(f_{j-1}\circ\ldots\circ f_i)=d_i-\dim_{\mathbb{K}}{\rm Hom}_Q(U_{i,j-1},V)&\\&=d_i-{\rm hom}(\alpha_{i,j-1},{\bf d})=d_i-\max\{d_i-d_k\, :\, i\leq k\leq j\}=\min\{d_i,\ldots,d_j\}.\end{align*}
The above formula thus constructs a representation $V$ where all the ranks $r_{ij}(V)$ reach their theoretical maximum $\min\{d_i,\ldots,d_j\}$.
\subsection{} Now we consider the case when $Q$ has a unique sink, that is, 
$$1\rightarrow 2\rightarrow\cdots\rightarrow s\leftarrow \cdots \leftarrow n-1 \leftarrow n.$$
We thus find with our main result that for all $1\leq i\leq j\leq n$ we have (ignoring the redundant elements not contributing to the minimum)
$$\mathbf{m}_{V_{{\bf d},\rm{ri}}}(\alpha_{ij})=\begin{cases}\big[\min\big\{d_k+d_{k'}-d_s, d_s-d_{i-1}-d_{j+1},\\
\ \ \  \ d_{k'}-d_{j+1},d_{k}-d_{i-1}: i\leq k\leq s\leq k'\leq j
\big\}\big]_+,&\ \ \text{$i<s$ and $j>s$}\\
\left[\min\left\{d_k-d_{j+1}, d_k-d_{i-1}: i\leq k\leq j\right\}\right]_+,&\ \ \text{$j<s$ or $i>s$}\\
\left[\min\left\{d_k-d_{i-1}: i\leq k\leq j\right\}\right]_+,&\ \ \text{$j=s$}\\
\left[\min\left\{d_k-d_{j+1}: i\leq k\leq j\right\}\right]_+,&\ \ \text{$i=s$}
\end{cases}$$
%$$\mathbf{m}_{V_{{\bf d},\rm{ex}}}(\alpha_{ij})=\big[\min\big\{d_k+d_{k'}-d_s, \ d_{\ell'}-d_{j+1}, d_{\ell}-d_{i-1},d_{\ell'}+d_{\ell}-d_{j+1}-d_{i-1},
%$$$$d_s-d_{i-1}-d_{j+1}: i\leq \ell<s, s<\ell'\leq j, i\leq k\leq s\leq k'\leq j
%\big\}\big]_+,\ \ \text{ if $i<s$ and $j>s$}$$
%$$\mathbf{m}_{V_{{\bf d},\rm{ex}}}(\alpha_{ij})=\left[\min\left\{d_k-d_{j+1}, d_k-d_{i-1}: i\leq k\leq j\right\}\right]_+\ \ \text{ if $j<s$ or $i>s$}$$
%$$\mathbf{m}_{V_{{\bf d},\rm{ex}}}(\alpha_{ij})=\left[\min\left\{d_k-d_{i-1}: i\leq k\leq j\right\}\right]_+\ \ \text{ if $j=s$}$$
%$$\mathbf{m}_{V_{{\bf d},\rm{ex}}}(\alpha_{ij})=\left[\min\left\{d_k-d_{j+1}, d_k-d_{i-1}: i\leq k\leq j\right\}\right]_+\ \ \text{ if $i>s$}$$
%$$\mathbf{m}_{V_{{\bf d},\rm{ex}}}(\alpha_{ij})=\left[\min\left\{d_k-d_{j+1}: i\leq k\leq j\right\}\right]_+\ \ \text{ if $i=s$}$$
where we understand again $d_0=d_{n+1}=0$.
%For example, if $\bf{d}=\mathbf{dim}(\mathbb{K}Q\oplus (\mathbb{K}Q^*))$ we simply have $V_{{\bf d},\rm{ex}}\cong I(r)^2\oplus P(1)^{s-1}\oplus P(n)^{n-s}$. 
The representation $V=V_{\bf d, \rm{ri}}$, given as 
$$V_1\stackrel{f_1}{\longrightarrow} V_2\stackrel{f_2}{\longrightarrow}\cdots\stackrel{f_{s-1}}{\longrightarrow} V_s\stackrel{f_s}{\longleftarrow} \cdots \stackrel{f_{n-2}}{\longleftarrow} V_{n-1} \stackrel{f_{n-1}}{\longleftarrow} V_n$$
is again determined up to isomorphism by the ranks $r_{ij}(V)$, $1\leq i\leq j\leq n$, where $r_{ii}(V)=\mathrm{dim}_{\mathbb{K}}(V_i)$ as above and otherwise 
$$r_{ij}(V)=\begin{cases}{\rm rank}(f_{j-1}\circ\cdots\circ f_i:V_i\rightarrow V_{j}),& \text{if $j\leq s$}\\
{\rm rank}(f_{i}\circ\cdots\circ f_{j-1}:V_j\rightarrow V_{i}),& \text{if $i\geq s$}\\
{\rm rank}((f_{s-1}\circ\cdots\circ f_i)\oplus (f_{s}\circ\cdots\circ f_{j-1}):V_i\oplus V_j\rightarrow V_r),& \text{if $j>s, i<s$}\\
%$$r_{ij}(V)=\begin{cases}{\rm rank}(f_{j}\circ\cdots\circ f_i:V_i\rightarrow V_{j+1}),& \text{if $j<s$}\\
%{\rm rank}(f_{i-1}\circ\cdots\circ f_{j-1}:V_j\rightarrow V_{i-1}),& \text{if $j>s, i>s$}\\
%{\rm rank}((f_{s-1}\circ\cdots\circ f_i)\oplus (f_{s}\circ\cdots\circ f_{j-1}):V_i\oplus V_j\rightarrow V_r),& \text{if $j>s, i<s$}\\
\end{cases} $$
We have 
$$r_{ij}(V)=\min\{d_i,\ldots,d_j\},\ \ \text{if $j\leq s$ or $i\geq s$}$$
as in the equioriented case, and for $i<s<j$ we get
\begin{align*}r_{ij}(V)&=d_i+d_j-\dim_{\mathbb{K}}{\rm Hom}_Q(U_{i,j},V)=d_i+d_j-{\rm hom}(\alpha_{i,j},{\bf d})&\\&
=d_i+d_j-\mathrm{max}\{d_i+d_j-d_k-d_{k'}, d_i+d_j-d_s: i\leq k\leq s, s\leq k'\leq j\}&\\&=\mathrm{min}\{d_k+d_{k'}, d_s: i\leq k\leq s, s\leq k'\leq j\}&\\&=\mathrm{min}\left\{\mathrm{min}\{d_i,\dots,d_s\}+\mathrm{min}\{d_s,\dots,d_j\}, d_s\right\}.\end{align*}
So in both cases, the rigid representation is the unique representation maximizing certain ranks. This will be generalized to arbitrary orientations in the next subsection and can be interpreted as a coherent linear algebra description of the rigid representation.

\subsection{} Consider an arbitrary orientation of the type $A$ quiver $Q$ and let $Q(i,j)$ for $i\leq j$ be the full subquiver of $Q$ with vertices $\{i,\dots,j\}$. Moreover, we introduce the following notation.
\begin{itemize}
\item Let $Q^{\rm{si}}(i,j)$ (resp. $Q^{\rm{so}}(i,j)$) be the set of sinks in $Q(i-1,j+1)$ different from the nodes $i$ and $j$ (resp. the set of sources in $Q(i,j)$).\vspace{0,1cm}

\item For $k\in Q^{\rm{si}}(i,j)$ we denote by $w_k^{ij}$ (resp. $u_k^{ij}$) the unique longest path in $Q(i,k)$ (resp. $Q(k,j)$) ending in $k$.
\end{itemize}
%denote by $Q^{\rm{si}}(i,j)$ (resp. $Q^{\rm{so}}(i,j)$) the set of sinks in $Q(i,j)$ different from $i+1$ and $j-1$ (resp. the set of sources in $Q(i,j)$). 
Then we have a projective resolution of $U_{i,j}$:
$$0\longrightarrow P_1:=\bigoplus_{k\in Q^{\rm{si}}(i,j)}P(k)\stackrel{\varphi}{\longrightarrow} P_0:=\bigoplus_{\ell\in Q^{\rm{so}}(i,j)}P(\ell)\longrightarrow U_{i,j}\longrightarrow 0$$
where $P(k)$ denotes as usual the indecomposable projective module corresponding to the node $k$. The map $\varphi$ restricted to $P(k)$ is given by $\varphi(e_{k})=w_k^{ij}-u_k^{ij}$ (where $e_k$ is the constant path at node $k$).
Now for any $V\in R_{\bf d}(Q)$ the functor $\mathrm{Hom}_Q(-,V)$ yields: 
$$0\rightarrow \mathrm{Hom}_Q(U_{i,j},V)\rightarrow \mathrm{Hom}_Q(P_0,V)\xrightarrow[]{\Phi^V_{i,j}} \mathrm{Hom}_Q(P_1,V)$$
and hence
\begin{equation}\label{glh42}\mathrm{dim}_{\mathbb{K}}\mathrm{Hom}_Q(U_{i,j},V)=\mathrm{dim}_{\mathbb{K}}\mathrm{Hom}_Q(P_0,V)-\mathrm{rank}(\Phi^V_{i,j}).\end{equation}
%$$H(U_{i,j},M)=H(P_0,M)-\mathrm{rk}(\Phi^M_{i,j}).$$
%In particular, we get for all $N,M\in R_{\bf d}(Q)$ (using \cite{}):
%$$\mathcal{O}_N\subseteq \overline{\mathcal{O}}_M\iff \mathrm{dim}_{\mathbb{K}}(\mathrm{Hom}_Q(U_{ij},M))\leq \mathrm{dim}_{\mathbb{K}}(\mathrm{Hom}_Q(U_{ij},N))\ \forall 1\leq i\leq j\leq n$$ 
%$$\iff \mathrm{rk}(\Phi^N_{i,j})\leq \mathrm{rk}(\Phi^M_{i,j}) \ \forall 1\leq i\leq j\leq n.$$
%$P_1=\bigoplus_{k\in Q^{\rm{si}}(i-1,j+1)}P(k)$ and $P_0=\bigoplus_{k\in Q^{\rm{so}}(i,j)}P(k)$ 

%where $\{k_1<\cdots <k_r\}\subseteq Q(i-1,j+1)$ is the set of sinks in $Q(i-1,j+1)$ different from $i$ and $j$ and $\{\ell_1<\cdots <\ell_p\}\subseteq Q(i,j)$ is the set of sources in $Q(i,j)$
%$P_1=P(k_1)\oplus\cdots \oplus P(k_r)$ and $P_0=P(\ell_1)\oplus\cdots \oplus P(\ell_p)$

%\begin{lemma}
%Let $\{k_1<\cdots <k_r\}\subseteq Q(i-1,j+1)$ the set of sinks in $Q(i-1,j+1)$ different from $i$ and $j$ and $\{\ell_1<\cdots <\ell_p\}\subseteq Q(i,j)$ the set of sources in $Q(i,j)$. Then we have a projective resolution 
%$$0\rightarrow P_1\rightarrow P_0\rightarrow U_{ij}\rightarrow 0$$
%where $P_1=P(k_1)\oplus\cdots \oplus P(k_r)$ and $P_0=P(\ell_1)\oplus\cdots \oplus P(\ell_p)$.
%\begin{proof}
%The map $\varphi:P_1\rightarrow P_0$ restricted to $P(k)$ is given by $\varphi(e_{k})=w_k^1-w_k^2$ ($e_k$ is the constant path at node $k$) where $w_k^1$ (resp. $w_k^2$) is the unique longest path in $Q(i,k)$ (resp. $Q(k,j)$) ending in $k$.
%; recall that $k$ is a sink with $k\neq i,j$.
%\end{proof}
%\end{lemma}
Recall the elements $w_k^{ij}$ and $u_k^{ij}$ from above and set $r_{ij}(V):=\mathrm{rank}(A_{ij}(V))$ for all $1\leq i\leq j\leq n$ where
\begin{equation}\label{rk130} A_{ij}(V):=\left[ \begin{array}{rrrrrr}
V_{w_{k_1}^{ij}} &V_{u_{k_1}^{ij}}& 0&\cdots &0 &0\\ 
0 & V_{w_{k_2}^{ij}} & V_{u_{k_2}^{ij}}&\cdots &0&0 \\
\vdots & \vdots &\vdots&  &\vdots&\vdots\\
0 & 0 & 0& \cdots &V_{w_{k_r}^{ij}} & V_{u_{k_r}^{ij}}\\ 
\end{array}\right]\end{equation}
and $Q^{\rm{si}}(i,j)=\{k_1<\cdots <k_r\}$.
\begin{proposition}
Let $V\in R_{\bf d}(Q)$. Then $V\cong V_{\bf d,\rm{ri}}$ if and only if 
\begin{equation}\label{rk13}r_{ij}(V)= \sum_{\ell \in Q^{\rm{so}}(i,j)}d_{\ell}-{\rm hom}(\alpha_{i,j},{\bf d}),\ \text{ for all $1\leq i\leq j\leq n$.}\end{equation}
\begin{proof}
We claim that $r_{ij}(V)= \mathrm{rank}(\Phi^V_{i,j})$ for all $1\leq i\leq j\leq n$. Before we prove the claim, let us see how this proves the proposition. From \eqref{glh42} we have  
$$\mathrm{rank}(\Phi^{V_{\bf d,\rm{ri}}}_{i,j})=\sum_{\ell \in Q^{\rm{so}}(i,j)}d_{\ell}-{\rm hom}(\alpha_{i,j},{\bf d}).$$ 
Conversely, if $V$ satisfies \eqref{rk13}, we get again with \eqref{glh42} that $\mathrm{rank}(\Phi^V_{i,j})= \mathrm{rank}(\Phi^{V_{\bf d, \rm{ri}}}_{i,j})$ for  $1\leq i\leq j\leq n$, thus
$$\mathrm{dim}_{\mathbb{K}}\mathrm{Hom}_Q(U_{i,j},V)=\mathrm{dim}_{\mathbb{K}}\mathrm{Hom}_Q(U_{i,j},V_{\bf d, \rm{ri}})\, , 1\leq i\leq j\leq n.$$
Hence $V\cong V_{\bf d,\rm{ri}}$ since representations are determined up to isomorphism by the dimensions of homomorphism spaces from indecomposables. So it remains to show $r_{ij}(V)= \mathrm{rank}(\Phi^V_{i,j})$.
The rank of $\Phi^V_{i,j}$ is obviously given by the dimension of the space $\{g\circ \varphi: g\in \mathrm{Hom}_Q(P_0,V)\}$. Note that any $g\in \mathrm{Hom}_Q(P(k),V)$ is determined by an element $V_k^g\in V_k$ via $g(e_k)=V_k^g$. Let $Q^{\mathrm{so}}(i,j)=\{\ell_1<\cdots <\ell_p\}$ and define a map 
$$\{g\circ \varphi: g\in \mathrm{Hom}_Q(P_0,V)\}\rightarrow  \mathrm{Im}\Big(\sum_{m=1}^r V_{w_{k_m}^{ij}}\oplus-V_{u_{k_m}^{ij}}: V_{\ell_1}\oplus\cdots \oplus V_{\ell_p}\rightarrow V_{k_1}\oplus \cdots \oplus V_{k_r}\Big)$$
as follows. Recall that $$V_{w_{k_m}^{ij}}\oplus-V_{u_{k_m}^{ij}}:V_{h(w^{ij}_{k_m})}\oplus V_{h(u^{ij}_{k_m})}\rightarrow V_{k_m}$$
and define the assignment 
$$g\circ \varphi\mapsto \sum_{m=1}^r V_{w_{k_m}^{ij}}(V^g_{h(w^{ij}_{k_m})})-V_{u_{k_m}^{ij}}(V^g_{h(u^{ij}_{k_m})})= \sum\limits_{m=1}^{r}g\circ \varphi (e_{k_m}).$$
Recall that the starting vertices of $w_k^{ij}$ and $u_k^{ij}$, respectively, correspond to sources, that is, the projective module satisfies $P(h(w_k^{ij}))\subseteq P_0$. This map is in fact a vector space isomorphism. To see this let $\beta\in P(k)\subseteq P_1$ be a path starting at $k$. Then $$(g\circ \varphi)(\beta)=V_{\beta}(g\circ\varphi(e_k))=V_{\beta}(g(w^{ij}_{k})-g(u_{k}^{ij}))=V_{\beta}(V_{w_{k}^{ij}}(V^g_{h(w^{ij}_{k})})-V_{u_{k}^{ij}}(V^g_{h(u^{ij}_{k})}))$$
and the injectivity is established. The surjectivity is clear, since given an element $$\sum_{m=1}^r V_{w_{k_m}^{ij}}(z_{h(w^{ij}_{k_m})})-V_{u_{k_m}^{ij}}(z_{h(u^{ij}_{k_m})})$$ in the image we can simply define $g\in \mathrm{Hom}_Q(P_0,V)$ by $g(e_{h(w_k^{ij})})=z_{h(w^{ij}_{k})}$ and $g(e_{h(u_k^{ij})})=z_{h(u^{ij}_{k})}$ respectively. The map \(\sum_{m=1}^r V_{w_{k_m}^{ij}}\oplus-V_{u_{k_m}^{ij}}\) is represented by the matrix \(A_{i,j}(V)\) up to some minus signs, which are negligible for the rank.
\end{proof}
\end{proposition}
So the rigid representation $V_{\bf d, \rm{ri}}$ is the unique representation up to isomorphism in $R_{\bf d}(Q)$ maximizing the rank tuple $(r_{ij}(V))_{1\leq i\leq j\leq n}$ coordinatewise.\\[1ex]
We conclude this section with a few remarks and further questions:

\begin{itemize}
\item The piecewise-linear nature of our formula is natural a priori, from the fact that the functions ${\bf m}_{V_{{\bf d},{\rm ri}}}$ have to be supported on one of finitely many maximal subset of roots with vanishing ${\rm Ext}^1_Q$-groups between the corresponding indecomposable representations, which feature on the Catalan combinatorial aspects of quivers \cite{RingelCatalan}. 
\item The minimization in Theorem \ref{main} includes many redundancies, and it would be interesting to find an optimal form of the formula.
\item Due to the Auslander-Reiten theoretic nature of the proof of Theorem \ref{main}, generalizations to the case of preprojective or preinjective dimension vectors for arbitrary acyclic quivers can be expected.
\end{itemize}

%\begin{definition}
%Given a representation $M\in R_{\bf d}(Q)$, we define its \textit{rank tuple} as follows: $$\mathbf{rk}(M)=(r_{i,j}(M))_{1\leq i\leq j\leq n},\ \ \ r_{i,j}(M)=\mathrm{rk}(A_{i,j}(M)).$$
%\end{definition}
%So the above formula constructs a representation $V$ where all ranks $r_{i,j}(V)$ reach their maximum.
%\begin{example}
%Let $Q$ be the $A$ quiver with a unique sink at node $r$ and $\bf{d}=\mathbf{dim}(\mathbb{K}Q\oplus (\mathbb{K}Q^*))$. Then the exceptional representation is given by $V_{\bf d, \rm{ex}}=I(r)^{\oplus 2}\oplus P(1)^{r-1}\oplus P(n)^{n-r}$. We will argue this using both, the general formula and the rank tuples. To see this we first note that 
%$${\bf d}=(r+1,\dots,r+1,n+1,n-r+2,\dots,n-r+2).$$
%Hence, to maximize the ranks we have to choose the identity map $f_{\alpha}=\mathrm{id}$ for all $\alpha\in Q_1$ except for the arrows at the unique sink
%$$r-1\xrightarrow{f_1} r\xleftarrow{f_2} r+1$$
%Both maps have to be embeddings and the rank of the matrix $(f_1 f_2)$ has to be maximized. So if $f_1$ is represented by the matrix
%$$\left[ \begin{array}{l}
%\mathds{1}_{r+1} \\ 
%\mathbb{O} \\
%\end{array}\right]$$
%than the maximal rank is obtained for the map $f_2$ represented by 
%$$\left[\begin{array}{l}
 %\mathbb{O}  \\ 
%\mathds{1}_{n-r+2}  \\

%\end{array}\right]$$
%......
%%\end{example}

%We apply this in particular to the dimension vector ${\bf d}={\rm\bf dim}(kQ\oplus kQ^*)\ldots$

\bibliography{openorbit}{}
\bibliographystyle{plain}

\end{document}